\documentclass{mcom-l}
\usepackage{array,amsxtra}
\usepackage{amsmath}
\usepackage{amsfonts}
\usepackage{amssymb}
\usepackage{amsthm}
\usepackage{bbm}
\usepackage{graphicx}
\usepackage{psfrag}
\usepackage{cite}
\usepackage{hyperref}
\usepackage{cleveref}
\hypersetup{
urlcolor=red}
\usepackage{url}
\usepackage[usenames]{color} 
\usepackage{epsfig}
\usepackage{mathrsfs} 

\usepackage{cleveref}
\usepackage{latexsym}
\usepackage{mathtools}
\usepackage{listings}
\newtheorem{theorem}{Theorem}[section]
\newtheorem{proposition}[theorem]{Proposition}

\newtheorem{lemma}[theorem]{Lemma}
\newtheorem{corollary}[theorem]{Corollary}
\newtheorem{remark}[theorem]{Remark}

\newcommand{\spn}{\mathop{\mathrm{span}}}

\newcommand{\nats}{\mathbb{N}}
\newcommand{\reals}{\mathbb{R}}
\newcommand{\RR}{\mathbb{R}}

\newcommand{\CC}{\mathbb{C}}

\newcommand{\calo}{\mathcal{O}}

\newcommand{\sph}{\mathbb{S}} 

\newcommand{\bfa}{\mathbf{a}}
\newcommand{\bfb}{\mathbf{b}}
\newcommand{\bfdel}{\boldsymbol{\delta}}

\newcommand{\bfi}{\mathbf{i}}
\newcommand{\bfj}{\mathbf{j}} 
\newcommand{\bfk}{\mathbf{k}}
\newcommand{\bfn}{\mathbf{n}}

\newcommand{\bfx}{\mathbf{x}}

\renewcommand{\d}{\mathrm{dist}}

\renewcommand{\L}{\mathcal{L}}

\newcommand{\dist}{\operatorname{dist}}

\newcommand{\cali}{\mathcal I}

\newcommand{\Omegal}{\overline \Omega}

\newcommand{\STRset}[1]{\{#1\}}

\numberwithin{equation}{section}

\title[A Meshless Kernel Method for Non-Local Diffusion]{A Meshless Galerkin Method for Non-Local Diffusion Using Localized Kernel Bases}
\author{R. B. Lehoucq}
\address{Computational Mathematics, Sandia National Laboratories,  Albuquerque, NM 87185-1320.} 
\email{rblehou@sandia.gov}
\thanks{Research supported by the Laboratory Directed Research and
Development (LDRD) program at Sandia National Laboratories. Sandia is
multi-program laboratory managed and operated by Sandia Corporation, 
wholly a subsidiary of Lockheed Martin Corporation, for the U.S. Department
of Energy's National Nuclear Security Administration under contract
DE-AC04-94AL85000.}

\author{F. J.~Narcowich}
\address{ Department of Mathematics, Texas A\&M University, College
  Station, TX 77843. } 
\email{fnarc@math.tamu.edu}
\thanks{Research supported by grant DMS-1514789 from the National
  Science Foundation.}

\author{S. T. Rowe}
\address{Sandia National Laboratories,  Albuquerque, NM 87185.} 
\email{srowe@sandia.gov}
\thanks{Research supported by grant DMS-1211566 from the National Science Foundation and 
Sandia National Laboratories.} 

\author{J. D. Ward}
\address{ Department of Mathematics, Texas A\&M University,
    College Station, TX 77843.} 
\email{jward@math.tamu.edu} 
\thanks{Research supported by grant DMS-1514789 from the National Science
  Foundation.}



\begin{document}
\maketitle

\begin{abstract}
We introduce a meshless method for solving both continuous and discrete variational formulations of a volume constrained, nonlocal diffusion problem. We use the discrete solution to approximate the continuous solution. Our method is nonconforming and uses a localized Lagrange basis that is constructed out of radial basis functions. By verifying that certain inf-sup conditions hold, we demonstrate that both the continuous and discrete problems are well-posed, and also present numerical and theoretical results for the convergence behavior of the  method. The stiffness matrix is assembled by a special quadrature routine unique to the localized basis. Combining the quadrature method with the localized basis produces a well-conditioned, symmetric matrix. This then is used to find the discretized solution.
\end{abstract}

\section{Introduction}
\label{intro}

The contribution of our paper is a rigorous numerical analysis of a meshless method for solving  a variational formulation of a volume constrained, nonlocal diffusion problem. Our method is nonconforming and uses a localized Lagrange basis that is constructed out of radial basis functions.  The analysis presented demonstrates that the Lagrange multiplier method introduced in \cite{bolr:15} for nonlocal diffusion  is well posed, in both the discrete and continuous cases. Our paper also replaces the Lagrange functions considered in \cite{bolr:15} with local Lagrange functions as in \cite{lero:15}, leading to dramatically reduced quadrature costs. 

Nonlocal diffusion generalizes classical diffusion by replacing the partial differential equations with integral equations.  Various models have been proposed for these
cases of so-called anomalous diffusion, which include models based on integral equations and fractional derivatives. The nonlocal equation we consider has applications in a
variety of fields besides anomalous diffusion such as image analyses, nonlocal heat conduction, machine learning, and peridynamic mechanics. We apply our radial basis method to a volume constrained diffusion equation. Volume constraints replace the boundary conditions associated with classical diffusion, and are needed to demonstrate that the problem is well posed. It also provides a link with a Markov jump process; see \cite{dglz:12,budl:14} for additional information, motivation and citations to the literature. 

An important distinction with the radial basis method introduced in \cite{lero:15} and that of this paper, is that the former method is conforming whereas the latter is nonconforming, an unavoidable aspect of a fully radial basis function method given a volume constraint. The nonconforming method of local Lagrange functions then enjoys all the benefits of a radial basis function method -- error estimates and stability estimates. This represents a powerful manner in which a class of radial basis function methods can be used to approximate the solution of conventional weak formulations of classical boundary value problems. 

Meshfree methods obviate the need to mesh the domain. As noted in \cite{babuaska_survey03}, the development of meshless methods was stimulated
by difficulties related to mesh generation such as when the underlying domain has a complicated geometry or when remeshing is required for time-dependent problems.
Also mentioned in \cite{babuaska_survey03} was the potential advantages of meshless methods when a Lagrangian formulation is employed, which will be the case for this paper. Meshless methods also allow for flexibility
in the selection of approximating functions, in particular non-polynomial approximating functions. In this paper the approximating spaces will be spanned by certain localized kernel bases \cite{FHNWW2012, HNRW_15} that are 
distinguished by a rigorous approximation theory and give rise to  very practical and efficient numerical methods.

A conforming discontinuous Galerkin method for a nonlocal diffusion problem was introduced in \cite{dglz:12} where the basis functions are given by discontinuous piecewise polynomials. Assembly of this stiffness matrix results in a challenging problem in quadrature for two reasons. The first is that there are iterated integrals over $2n$ dimensional regions, where $n$ is underlying the spatial dimension, and the second is that the regions of integration involve partial element volumes. In contrast, the primary advantage of the meshfree methods is that entries in the stiffness matrix only require a pointwise evaluation of the kernel and multiplication by quadrature weights--complications arising from overlapping partial element volumes are irrelevant. Consequently, our proposed method requires only information at the radial basis function nodes or centers and also yields a straight forward assembly of a sparse stiffness matrix.  

The numerical analysis provided in this paper will be based on two specific classes of local Lagrange functions that will play the role of bases for the spaces $U_h$ and $\Lambda_h$ appearing in \eqref{bases}.
In \cite{HNRW_15}, it was shown  that for either thin-plate splines or Mat\'ern kernels on $\mathbb{R}^n$,  local Lagrange functions with each function determined by ${\calo}(\log N)^n$ points contained in a ball of radius $Kh\log h$ centered at a given point $\xi$ have very rapid decay around $\xi$. Moreover such functions generate very stable bases.

The theoretical development for such functions first appeared in \cite{FHNWW2012} in the context of $\sph^2$. The corresponding theory for compact domains in $\RR^n$ appeared in \cite{HNRW_15}.  Applications using these basis functions in the context of numerical solution of certain PDEs have been given in \cite{NRW_14,bolr:15,lero:15}. In particular stability estimates for this class of functions will play a crucial role in \Cref{discrete_case} for the numerical solvability of our problem.

The remainder of the paper is organized as follows. In section~\ref{var-form}, the variational framework for both the continuous and discrete cases is discussed; in addition, notation to be used throughout the paper is introduced. Section~\ref{rbfs} contains a review of the radial basis functions (RBFs) that give rise to the local Lagrange bases mentioned earlier. These bases are highly localized and computationally inexpensive. The main result is Theorem~\ref{q-interp_approx}, which provides Sobolev error estimates when approximation by the quasi-interpolation operator associated with the local Lagrange basis. Section~\ref{coercivitiy} establishes coercivity results for the bilinear form \eqref{bilinear-form}. 

The main results of the paper are presented in section~\ref{lag_mult_form}. The solutions to the Euler-Lagrange formulation \eqref{variationalform}, for both the continuous and discrete cases, are given in Theorem~\ref{lagrange_dirichlet_solution} and Theorem~\ref{lagrange_form_h_prop}, respectively. Finally, in section~\ref{numerical_results} numerical results are presented. These  results are in good agreement with the theoretical results discussed in section~\ref{sec:Error_Estimates}.

\section{Variational Formulation}
\label{var-form}

Consider a domain $\Omegal=\Omega\cup\Omega_{\mathcal{I}}$, where $\Omega$ is an inner domain, $\Omega_{\mathcal{I}}$ is the interaction region, and then define the bilinear form
\begin{equation} \label{bilinear-form}
a(u,v) = \int_{\Omegal}  \int_{\Omegal} \gamma(x,y)\big(u(x)-u(y)\big)\big(v(x)-v(y)\big)\, dx\,dy,
\end{equation}
where $\gamma(x,y)\ge 0$ is in $L_\infty(\Omegal,\Omegal)$ and $u,v\in L_2(\Omegal)$.  We assume that there exists an $L_\infty$ function  $\gamma_\delta :[0,\infty)\to [0,\infty)$, with support in $0\le r \le \delta<\infty$,  and that there are constants $c_1$ and $c_2$ for which
\begin{equation}
\label{anisotropic_bnds}
c_1\gamma_\delta(|x-y|) \le \gamma(x,y) \le c_2 \gamma_\delta(|x-y|), \ x,y\in \Omegal.
\end{equation}
Suppose also that $f\in L_2^c(\Omegal)$ where 
\begin{equation}
L_2^c(\Omegal) = \{ f \,|\, f \in L_2(\Omegal) \text{ with } f|_{\Omega_{\mathcal{I}}}=0 \}\,.
\end{equation}
Denote the inner product and norm on $L_2(\Omegal)$ by $\langle \cdot ,\cdot \rangle_{\Omegal}$ and $\| \cdot \|_{\Omegal}$, respectively. We will use similar notation for $L_2(\Omega)$ and $L_2(\Omega_\cali)$. 

We define the energy functional $E$ by
\begin{equation} \label{energy}
\left\{
\begin{aligned}
E(u) = \frac{1}{2} a(u,u) - \langle u ,f \rangle_{\Omegal}, \\
\text{subject to } u = 0 \text{ over } \Omega_{\mathcal{I}}\,.
\end{aligned}
\right.
\end{equation}
The  constraint over the volume $ \Omega_{\mathcal{I}}$ is the nonlocal analogue of a homogenous Dirichlet boundary condition; the reader is referred to \cite[pp.678--680]{dglz:12} for details and discussion. The paper \cite{dglz:12} demonstrated that the problem of finding the minimum of the energy functional was shown to be well-posed for $u$ in an energy constrained space $L^2_c(\Omegal) \subsetneq L_2(\Omegal)$.  In contrast, as in \cite{bolr:15}, we minimize the functional by the method of Lagrange multipliers because the local Lagrange basis is not contained in the energy constrained space.  The Lagrangian is defined as
\begin{align*}
L(u,\lambda) = E(u) +  b(u,\lambda), \ \text{where } b(u,\lambda):=\langle u , \lambda \rangle_{\Omega_{\mathcal{I}}}.
\end{align*}
Here, $\lambda \in L^2(\Omega_{\mathcal{I}})$ is the Lagrange multiplier.


The Euler-Lagrange formulation of the problem is then: Find $u \in L^2(\Omegal )$ such that 
\begin{equation}
\left\{\begin{aligned}
a(u,v)   + \langle u , \lambda \rangle_{\Omega_{\mathcal{I}}}  &=\langle u ,f \rangle_{\Omega}  &\text{ for all } v \in L^2(\Omegal ) \,, \\
 \langle u , w \rangle_{\Omega_{\mathcal{I}}}  &= 0 &\text{ for all } w \in L^2(\Omega_\mathcal{I})\,.
\end{aligned} \right.
\label{variationalform}
\end{equation}

We discretize this system by choosing finite dimensional subspaces $U_h \subset L^2(\bar \Omega)$ and $\Lambda_h \subset L^2(\Omega_{\mathcal{I}})$ where 
\begin{equation}\label{bases}
U_h = \text{span}\{\phi_i\}_{i=1}^N \,, \qquad \Lambda_h = \text{span} \{\psi_k\}_{k=1}^{N_{\mathcal{I}}}\,.
\end{equation}
We then approximate the pair $(u,\lambda)$ by the discrete pair $(u_h,\lambda_h)$ given by the expansions
\[
u_h = \sum_{j=1}^N \alpha_j \phi_j\,, \qquad \lambda_h = \sum_{k=1}^{N_{\mathcal{I}}} \beta_k \psi_{k}\,.
\]
Inserting the expansions into \eqref{variationalform} and in turn selecting $v$ and $w$ equal to each $\phi_i $ and $\psi_{k}$, respectively, determines the needed coefficients as the solution to the saddle point system
\begin{subequations} \label{discrete-saddlepoint}
\begin{align}
\begin{pmatrix}
A & B \\ B^T & 0
\end{pmatrix}
\begin{pmatrix}
\alpha \\ \beta
\end{pmatrix}
=
\begin{pmatrix}
b \\ 0
\end{pmatrix},
\end{align}
with matrix, vector entries given by
\begin{gather}
A_{i,j} = a(\phi_i\,,\phi_j) \,, \qquad B_{i,k} =  \langle  \phi_i \,, \psi_{k} \rangle_{\Omega_{\mathcal{I}}}  \,, \qquad b_i = \langle \phi_i ,f \rangle_{\Omega}\,.
\end{gather}
\end{subequations}

\section{Radial Basis Functions and Localized Kernel Bases}
\label{rbfs}

In this section, we give background material on interpolation and approximation with radial basis functions (RBFs). Radial basis functions  are used to construct the approximation space for the Galerkin method we propose in section~\ref{discrete_case}. The interested reader should consult \cite{wendland_book} or \cite{Fasshauer-2007-1} for further details on radial basis functions and interpolation.

\subsection{Radial basis functions}\label{sec:RBF}
A radial basis function (RBF) is a radial function $\Phi(x) = \phi(|x|)$, where $\phi \in C[0,\infty)$, that is (strictly) positive definite on $\RR^n$ \cite[Chapter~6]{wendland_book} or (strictly) conditionally positive definite on $\RR^n$, with respect to the set of polynomials $\pi_{m-1}:=\pi_{m-1}(\RR^n)$ having total degree $m-1$ or less \cite[Chapter~8]{wendland_book}. Specifically, this means that for every distinct set 
$X \subset \reals^n$, with cardinality $|X|=N<\infty$, and all nonzero $a\in \RR^{N}$ satisfying $\sum_{\xi\in X} a_{\xi} p(\xi) = 0$,  we have that 
\[
\sum_{\xi \in X}\sum_{\zeta\in X} \phi(|\xi-\zeta|) a_{\xi} a_{\zeta} >0.
\]

The RBFs that are conditionally positive definite with respect to $\pi_{m-1}$ are said to have order $m\ge 1$. If an RBF is positive definite, it has order 0. 

We will be especially interested in \emph{thin-plate splines} (TPS) or \emph{surface splines}, because they produce Lagrange and local Lagrange functions that are well-localized in space and have a ``small'' footprint among the thin-plate splines used to construct them; see \cite{HNRW_15}. For $m>n/2$, a thin-plate spline  is defined by
\begin{equation}\label{surface_spline}
\phi_m(r):=
\begin{cases} r^{2m-n}&n \text{ is odd}\\
r^{2m-n}\log r&n \text{ is even}.
\end{cases}
\end{equation}
For each $m>n/2$, the TPS $\phi_m(|x|)$  is an order $m$ RBF.

An example of an order 0 RBF that has properties similar to a TPS is the Mat{\' e}rn kernel, which is defined by 
\begin{equation}
\label{Matern}
\kappa_m(r):= CK_{m-n/2}(r)\,r^{m-n/2},\ m>n/2.
\end{equation}
Here $C$ is a constant depending on $\mu$ and $n$, and $K_{\nu}$ is a Bessel function of the second kind. 

The approximation space for any RBF $\Phi(x)=\phi(|x|)$ of order $m$ associated with a unisolvent\footnote{Unisolvent with respect to $\pi_{m-1}$ means every $p\in \pi_{m-1}$  is uniquely determined by its values on $X$.} set $X$, which is called the set of centers, is defined by 
\begin{equation}
\label{approx_space}
V_X := \bigg\{\sum_{\xi \in X} a_\xi \phi(|x- \xi|) \colon \sum_{\xi\in X}a_\xi\,p(\xi) =
0 \ \forall\ p\in \pi_{m-1} \ \bigg\}+\pi_{m-1}, \ \text{where }\pi_{-1}=\{0\}. 
\end{equation}
Specifically, each $s\in V_X$ has the form\footnote{Bases other than $\{x^\gamma\}_{|\gamma| \le m-1}$ may be used for $\pi_{m-1}$.}
\begin{equation}
\label{RBF_interpolant}
s(x) = \sum_{\xi \in X}  a_\xi\phi(|x - \xi |) +\sum_{|\gamma| \le m-1} \beta_\gamma x^\gamma,
\end{equation}
where $\gamma=(\gamma_1,\ldots,\gamma_n)$ is a multi-index, $|\gamma|=\gamma_1+\cdots+\gamma_n$, and $\sum_{\xi\in X}a_\xi p(\xi)=0$ for all $p\in \pi_{m-1}$.

If $X$ is a unisolvent set for $\pi_{m-1}$ and $d_\eta \in \CC$ is given at each $\eta \in X$,  there is a unique $s\in V_X$ 
that interpolates the $d_\eta$'s -- i.e., $s(\eta)=d_\eta$. The coefficients for $s$ in \eqref{RBF_interpolant} are obtained by solving the $N+\dim \pi_{m-1}$ equations
\begin{equation}
\label{interp_system}
\left\{
\begin{gathered}
\sum_{\xi \in X}  a_\xi\phi(|\eta - \xi |) + \sum_{|\gamma|\le m-1}\beta_\gamma \eta^\gamma = d_\eta, \ \eta \in X \\
\sum_{\xi\in X}a_\xi \xi^\nu = 0 \ \forall \ |\nu|\le m-1.
\end{gathered}\right.
\end{equation}

If the data are generated by a continuous function $f$, then we use $I_Xf$ instead of $s$. Finally, if $p$ is a polynomial in $\pi_{m-1}$ and $d_\eta = p(\eta)$,  then $I_Xp= p$. Thus, interpolation from $V_X$ reproduces polynomials in $\pi_{m-1}$. 

\subsubsection{Geometry of the set of centers}
\label{geometry}

The geometry of the centers is important for estimating the approximation quality of the RBF interpolant and for estimating the condition number of the interpolation matrix. RBF interpolation offers the advantage of not requiring regular distributions of points; arbitrarily scattered centers produce invertible interpolation matrices for positive definite functions.   

Let $D$ be a bounded, Lipschitz domain\footnote{To avoid notation confusion, we use $D$ rather than $\Omega$, which is is standard.} and let $X \subset D \subset \mathbb{R}^n$ be a set of scattered centers. We define the \emph{fill distance} (or \emph{mesh norm}) $h$, the \emph{separation radius} $q$ and the \emph{mesh ratio} $\rho$ to be:
\begin{equation} \label{geom_X}
 h:=\sup_{x\in D} \d(x,X), \quad    q:=\frac12 \inf_{\xi \in X} 
\d(\xi, X\!\setminus \!\{\xi\}) , \quad  \rho:=\frac{h}{q}.
\end{equation}  
The mesh norm $h$ is the radius of the largest ball in $D$ that does not contain any centers. The separation radius $q$ is the radius of the largest ball that can be placed at a center without including any other center; it is thus half of the minimal pairwise distance between the centers. Finally, we define the mesh ratio to be $h/q$. Obviously, $\rho\ge 1$. 

The mesh ratio measures  the uniformity of the  distribution of $X$ in $D$. The larger $\rho$ is, the less uniform the distribution is. If $\rho$ is ``small'', then we say that the point set $X$ is quasi-uniformly distributed, or simply that $X$ is quasi-uniform. Geometrically, $\rho$ controls how the centers cluster as the number of points increases. We note that for the quasi-uniformly distributed collections of centers $\STRset{X_{h,q}}$, which we will encounter later, we do not require that any of the point sets are nested in another. 

Earlier we mentioned that for a unique interpolant from $V_X$ to exist, it was necessary that $X$ be unisolvent with respect to $\pi_{m-1}$.  For a Lipschitz domain, there is a constant $r_{m,\partial D}$ such that if $h\le r_{m,\,\partial D}$ then $X$ will be unisolvent \cite[Proposition~3.5]{Narcowich-etal-2005-1} -- i.e., unisolvency holds if $h$ is small enough. 

\subsubsection{Approximation power}
\label{approx_power}

RBF interpolation and approximation provide excellent approximation power when $X$ is quasi-uniformly distributed in  $D$. For RBFs with Fourier transforms that behave like 
\begin{equation}\label{Phi-FT}
c_1(1+\|\omega\|_2^2)^{-\tau} \le \widehat{\Phi}(\omega)\le 
c_2(1+\|\omega\|_2^2)^{-\tau},\qquad \omega\in\RR^n,
\end{equation}
or have a generalized Fourier transform that satisfies
\begin{equation}\label{Phi-FT2}
c_1\|\omega\|_2^{-2\tau} \le \widehat{\Phi}(\omega)\le 
c_2\|\omega\|_2^{-2\tau},\qquad \tau\in 2\mathbb N,\; \omega\in\RR^n\setminus\{0\},
\end{equation}
where we take $\tau>n/2$, we have the approximate rates in the result below.

\begin{theorem}[\!\!{\cite[Theorem~4.2]{Narcowich-etal-06-3}}]
\label{approx_rates} 
Suppose that $\phi$ is an RBF that satisfies \eqref{Phi-FT} or \eqref{Phi-FT2} and that $X$ is quasi-uniformly distributed in $D$, with separation radius $q$ and mesh norm $h$. If $\tau\ge\beta$, $\beta=k+s$ with $0\le s<1$ and $k>n/2$, and if $f\in W^\beta_2(D)$, then
\[
\|f-I_X f\|_{W_2^\mu(D)} \le C h^{\beta-\mu}
\rho^{\tau-\mu} \|f\|_{W_2^\beta(D)}, \qquad
0\le \mu\le \beta,
\]
where $I_Xf$ is given by \eqref{RBF_interpolant}.
\end{theorem}

The thin-plate splines satisfy \eqref{Phi-FT2}, and both Mat\'ern kernels and Wendland functions satisfy \eqref{Phi-FT}. (See \cite[Sections 8.3 \& 9.4]{wendland_book}.)

The motivation for Theorem~\ref{approx_rates} above was getting error rates for interpolation in spaces having functions ``rougher'' than the RBFs used, which are in the reproducing kernel Hilbert space or \emph{native space} \cite{wendland_book}. In cases where the functions are twice as smooth as those in the native space, Schaback \cite[Theorem~5.1]{Schaback-00-1} has shown that, under certain additional conditions, the approximate rate doubles. A result specifically for thin-plate splines is given in \cite[Theorem~11.36]{wendland_book}. On a sphere, where there is no boundary, rates can be doubled, with there being no additional conditions \cite[Section~A.1]{NRW_14}.

An earlier version of the theorem that applied only to thin-plate splines for integer cases was proved by Brownlee and Light \cite[Theorem~3.5]{Brownlee-Light-03}.

\subsubsection{Lagrange functions}
\label{full_local_lagrange}

Suppose that $\phi(|x|)$ is an order $m$ RBF and $X$ is a set of centers in $\RR^n$.  By \Cref{sec:RBF} above, we can find a unique interpolant $\chi_\xi(x)\in V_X$ such that $\chi_\xi(\eta)=\delta_{\xi,\eta}$. In words, $\chi_\xi$ is $1$ at $\xi$ and $0$ at the remaining points in $X$. The function  $\chi_\xi$ is called a \emph{Lagrange function} centered at $\xi$ and is given by
\begin{equation}
\label{lagrange_function}
\chi_\xi(x) = \sum_{\eta\in X} \alpha_{\eta,\xi} \phi(|x-\eta | ) + \sum_{|\gamma|\le m-1} \beta_{\gamma, \xi} x^\gamma,
\end{equation} 
 where the coefficients are determined by \eqref{interp_system}, with $d_\eta=\delta_{\xi,\eta}$. It is easy to see that the Lagrange functions $\{\chi_\xi\}_{\xi\in X}$ form a basis for $V_X$, and that every $s\in V_X$ may be uniquely expressed as 
\begin{equation}
\label{s_in_V}
s(x)=\sum_{\xi \in X} s(\xi)\chi_\xi(x).
\end{equation}

At this point, we specialize to the thin-plate splines $\phi_m$, given in \eqref{surface_spline} and the Matr\'ern kernels, defined in \eqref{Matern}. For these RBFs, there are three more important features of the Lagrange basis. The first is a decay property of the Lagrange functions. This is what makes them local. By \cite[eqn.~(3.6)]{HNRW_15}, if $x\in D$, then there exist positive constants\footnote{We have used $\mu=2\nu$ in \cite[eqn.~(3.6)]{HNRW_15}; $\nu$ is defined in \cite[eqn.~(3.7)]{HNRW_15}.} $\nu = \nu(\phi,n)$ and $C=C(\phi,n)$ such that
\begin{equation}
\label{ptwise}
|\chi_{\xi}(x)| \le C \rho^{m-d/2} \exp\bigg(-2\nu \frac{\min(\dist(x,\xi),\dist(\xi,\partial D)}{h}\bigg).
\end{equation}

The second is that, by \cite[eqn.~(3.7)]{HNRW_15}, the $\alpha_{\eta,\xi}$'s  in \eqref{lagrange_function} decay exponentially in the distance from $\eta$ to $\xi$:
\begin{equation}\label{coeff_decay}
|\alpha_{\eta,\xi}| \le 
C q^{d-2m} \mathrm{exp}\left(-{\nu} \frac{\dist(\eta,\xi)}{h}\right).
\end{equation}
Because of this decay, the $\chi_\xi$'s, which are given in \eqref{lagrange_function}, require only a relatively small number of the $\phi(\cdot - \eta)$'s to approximately calculate them. That is, the $\chi_\xi$'s have a small ``footprint'' in the space of kernels. In \cite[Section~7]{FHNWW2012}, similar decay in Lagrange functions constructed using spherical basis functions was used to construct a preconditioner for solving the equations for the $\alpha_{\eta,\xi}$'s. 

The third concerns stability of the Lagrange basis. We begin by defining the \emph{synthesis} operator $T: \CC^{|X |}\to V_X$ by $T\bfa = \sum_{\xi\in\Xi} a_{\xi} \chi_{\xi}=:s$. In other words, $T$ takes a set of coefficients $\{a_\xi \}_{\xi \in \Xi}$ and outputs a function $s\in V_\Xi$ satisfying $s(\xi) = a_\xi$. If we use the $\ell_p(X)$ norm for $\CC^{|X |}$ and $L_p(D)$ for $V_X$, then the stability of the basis, relative to these norms, is measured by comparing $\|\bfa\|_{\ell_p(X)}$ and $\|s\|_{L_p(D)}$, which we now do.

\begin{proposition}[\!\!{\cite[Eqn.~(3.3) \& Theorem~3.10]{HNSW_2_2011}, \cite[Proposition~3.2]{HNRW_15}}]
\label{full_lowercomparison}
Suppose $D\subset\RR^n$ is a bounded Lipschitz domain and let $\rho\ge 1$ be a fixed mesh ratio.
If the RBF is either a thin-plate spline or a Mat\'ern kernel, then there exist constants $c>0$ and $q_0>0$, depending  on $\rho $, so that for 
$X \subset D$ satisfying
$q<q_0$, $h/q\le \rho$, and $1\le p \le \infty$, we have
\begin{equation}\label{full_riesz_ineq}
c \left\| \bfa  \right\|_{\ell_p(X)}     
\le  
q^{-n/p} \| \textstyle{\sum_{\xi\in X}}\,a_{\xi}\chi_{\xi} \|_{L_p(D)} \le C \rho^{m+n/p} \left\|\bfa \right\|_{\ell_p(X)}.
\end{equation}
\end{proposition}

In the proposition above, the Lagrange functions are constructed using only centers in $D$. This isn't sufficient for the applications that we have in mind. What will need to work with is a set of Lagrange functions $\{\chi_\xi\}_{\xi\in X}$ constructed  from a \emph{larger} set of centers, $\widetilde X \supset X$: $\chi_\xi$ will be given by \eqref{lagrange_function}, but with the sum being over $\eta \in \widetilde X$ instead of $\eta \in X$. The result for the $\widetilde X$ case corresponding to the one for $X$ can be found in \cite[Proposition~3.2]{HNRW_15}. We will discuss properties of these Lagrange functions below, in connection with \emph{local} Lagrange functions.

\subsection{Local Lagrange functions and quasi interpolants}
\label{loc_lag_quasi_interp}
Finding the Lagrange functions $\{\chi_\xi\}_{\xi\in X}$ requires solving an $N\times N$ system of equations, where $N=|X|$, to obtain the  $\alpha_{\eta,\xi}$'s in \eqref{lagrange_function}. If $N$ is large, say 30,000, then finding the $\chi_\xi$'s essentially requires solving a $30,000\times 30,000$ system. This is a formidable task. It can however be gotten around by using a basis of \emph{local} Lagrange functions, which provide the same features as the Lagrange basis. The local Lagrange functions are obtained by solving $N$ relatively small linear systems.

\subsubsection{Local Lagrange functions}
\label{sec:local_lagrange_functions}

Local Lagrange functions were first introduced for use on the sphere \cite{FHNWW2012}, where decay properties and quasi-interpolation convergence rates were studied. The local Lagrange basis can be constructed in parallel by solving small (relative to the number of centers) linear systems.  

Recent work \cite{HNRW_15} has extended theoretical properties of the local Lagrange basis to compact domains in $\mathbb{R}^n$. They are constructed in the following way. 

We begin with a Lipschitz domain $D$ and a set of quasi-uniform centers $X$ in $D$; $X$ has mesh norm $h$, separation radius $q$, and mesh ratio $\rho$.  About each center $\xi\in X$, consider a ball $B_{\xi,r_h}$ centered at $\xi$ and having radius $r_h := K h |\log h|$, where $K>0$ is a parameter at our disposal. We also consider an enlarged region $\widetilde D = \{x\in \RR^n\colon \d(x,D) \le r_h\}$. This region is $D$ together with points in a band of width $r_h$. 

The reason for augmenting $D$ is to have a set that contains all of the $B_{\xi,r_h}$'s, so that problems with points near $\partial D$ can be ameliorated. Augmentation is unnecessary on a sphere, since it doesn't have a boundary.

The next step is to add centers $X_\text{extra}$ to the band $\widetilde D\setminus \! D$ in such a way that the mesh norm and separation radius of the new set $\widetilde X = X\cup X_\text{extra}$ are $h$ and $q$, or differ only slightly from $h$ and $q$ \cite[Section~2.3]{HNRW_15}. For each $\xi$, we use $\widetilde X$ to obtain a set of centers $\Upsilon_\xi :=B_{\xi,r_h}\cap \widetilde X$. In addition to the $\Upsilon_\xi$'s,  we define the quantity
\begin{equation}
\label{J_def}
 J=K\nu/2 +2n -4m-1,
\end{equation}
which will appear in the sequel in various error estimates.

The final step is to construct the ``local'' Lagrange function for $\Upsilon_\xi$, which we will define to be $b_\xi$; it has the form
\begin{equation}
\label{local_lagrange}
b_\xi(x) := \sum_{\eta\in \Upsilon_\xi} \alpha_{\eta,\xi}\phi(|x-\eta|) + \sum_{|\gamma| \le m-1} \beta_{\gamma,\xi\,}x^\gamma.
\end{equation}
where the coefficients are determined by the equations
\begin{equation}
\label{loc_lag_eqns}
b_\xi(\eta) = \delta_{\xi,\eta} \ \ \text{and } \sum_{\eta\in \Upsilon_\xi} \alpha_{\eta,\xi} p(\eta) = 0, \ \forall \ p\in \pi_{m-1}.
\end{equation}
This is the same form as that for $\chi_\xi$ in \eqref{lagrange_function}, except that only centers in $\Upsilon_\xi$ are used to construct $b_\xi$ and, in addition, that some of those centers may be from outside of $X$. 

The $b_\xi$'s are constructed using centers in $\widetilde X$, and it follows that they are in $V_{\widetilde X}$. However, they don't form a basis for $V_{\widetilde X}$. There are only $|X|$ of them. Since $\dim V_{\widetilde X}=|\widetilde X|>|X|$, there are too few of them to form a basis. Moreover, since some of the centers come from outside of $D$, not all of them are in $V_X$. Thus, we define a new space, $\widetilde V_X =  \spn\{ b_{\xi}: \xi \in X\}$, for which $\{b_\xi\}_{\xi \in X}$ \emph{is} a basis. 

As we mentioned at the end of the previous section, we will work with Lagrange functions -- the  $\chi_\xi$'s -- having $\xi\in X$, but constructed from centers in $\widetilde X$. We remark that $\chi_\xi \ne b_\xi$. However, they are close -- a fact that will prove important in the sequel.

The properties of RBFs guarantee that the coefficients in \eqref{local_lagrange} always can be solved for using the equations in \eqref{loc_lag_eqns}. Letting $N_\xi=|\Upsilon_\xi|$ and $\widetilde N_\xi =N_\xi + \dim(\pi_{m-1})$, we see that the system  is  $\widetilde N_\xi\times \widetilde N_\xi$. Estimating $N_\xi$ may be done by comparing volumes of $B_{\xi,r_h}$ and of $B_{\xi,q}$, which has only a single point $\xi$ in it. The result is
 \[
N_\xi \approx \text{vol}(B_{\xi,r_h})/\text{vol}(B_{\xi,q}) = r_h^n/q^n = K^n\rho^n |\log h |^n.
\]
The same comparison yields $N\approx \text{vol}(D)/\text{vol}(B_{\xi,q}) \approx C\rho^n h^{-n}$, equivalently, $h \approx N^{-1/n}$. It follows that $N_\xi \approx C (\log N)^n$. Since there are $N$ centers in $X$, determining all of the $b_\xi$'s requires solving $N$ systems that have approximately $(\log N)^n$ variables each, if the small number of $\beta$ variables are ignored. 

These systems are symmetric and can be solved in parallel.  Contrast this with determining the $\chi_\xi$'s. Doing that task requires solving $N$ equations having $N$ variables each. 

Concerning stability,  local Lagrange bases enjoy the same properties as the full Lagrange bases. The pertinent result, which is the analogue of Proposition~\ref{full_lowercomparison} for the local Lagrange case, is given below.

\begin{proposition}[\!\!{\cite[Proposition~4.12]{HNRW_15}}]
\label{local_lowercomparison}
Suppose $D\subset\RR^n$ is a bounded Lipschitz domain and let $\rho\ge 1$ be a fixed mesh ratio.
If the RBF is either a thin-plate spline or a Mat\'ern kernel, then there exist constants $c>0$ and $q_0>0$, depending  on $\rho $, so that for 
$X \subset D$ satisfying
$q<q_0$, $h/q\le \rho$, and $1\le p \le \infty$, we have
\begin{equation}\label{local_riesz_ineq}
c \left\| \bfa  \right\|_{\ell_p(X)}     
\le  
q^{-n/p} \| \textstyle{\sum_{\xi\in X}}\,a_{\xi}b_{\xi} \|_{L_p(D)} \le C \rho^{m+n/p} \left\|\bfa \right\|_{\ell_p(X)}.
\end{equation}
\end{proposition}

There is another result that we will need in the sequel. It involves an inequality established in the course of proving \cite[Theorem~4.11]{HNRW_15}.

\begin{lemma} \label{norm_diff_loc_full}
Suppose that $J>n$. Under the same assumptions made in Proposition~\ref{local_lowercomparison}, there is a constant $C>0$ 
such that,  for  $0\le \sigma \le m-(n/2-n/p)_+$ 
when $1\le p< \infty$ (or $\sigma \in \nats$ and  $0\le \sigma<m - n/2$ when $p=\infty$), the following holds:
\begin{equation}
\label{norm_diff_loc_full_ineq}
\| \textstyle{\sum_{\xi\in X}}\,a_{\xi} (b_{\xi} - \chi_{\xi})\|_{W_p^{\sigma}(D)} 
\le 
C h^{J-n(\frac{p-1}{p})}
\|\bfa\|_{\ell_p(X)}. 
\end{equation}
\end{lemma}
\begin{proof}
The inequality was established in the proof of \cite[Theorem~4.11]{HNRW_15}; it is \cite[eqn.~(4.11)]{HNRW_15}. 
\end{proof}

\subsubsection{Quasi interpolants and quadrature} 
Given a continuous function $f$ defined on $D$, we can construct a quasi interpolant for $f$ using the local Lagrange functions: $\widetilde I_X\!f = \sum_{\xi \in X}f(\xi)b_\xi$. Because the $b_\xi$'s are not \emph{full} Lagrange functions, when we evaluate $\widetilde I_Xf $ at $x=\eta \in X$, we only get $\widetilde I_X\! f(\eta) = \sum_{\xi \in X}f(\xi)b_\xi(\eta)$. If  $\eta \in \Upsilon_\xi$, then $b_\xi = \delta_{\xi,\eta}$. However, for $\xi'\ne \xi$, $b_{\xi'}(\eta)\ne \delta_{\xi',\eta}$. Consequently,  $\widetilde I_X\! f(\eta)$ will, in general, \emph{not} be equal to $f(\eta)$. The point is $\widetilde I_X\! f(\eta) \approx f(\eta)$. The following result extends Theorem~\ref{approx_rates} to the quasi-interpolant case, with centers outside of $D$.


\begin{theorem}
\label{q-interp_approx}
Let $k\in \nats$, $0\le s <1$, and $n/2<\beta = k+s\le m$. Suppose that  $\phi$ be a thin-plate spline $\phi_m$ or a Mat\'ern kernel $\kappa_m$. If $f\in W^\beta_2(D)$ is compactly supported in $D$ and $0\le \mu \le \beta$, then there is an $h_0$ and a sufficiently large $K$ such that for all $h\le h_0$ we have
\begin{equation}
\label{q-interp_approx_rate}
\|f- \widetilde I_X\!f \|_{W^\mu_2(D)} \le C h^{\beta-\mu}\|f \|_{W^\beta_2(D)}.
\end{equation}
\end{theorem}

\begin{proof}
Let $x\in D$ and form both the interpolant $I_Xf = \sum_{\xi \in X} f(\xi)\chi_\xi(x)$ and the quasi interpolant $\widetilde I_X\! f(x) = \sum_{\xi \in X} f(\xi)b_\xi(x)$ for $f$. We have
\[
\|f- \widetilde I_X\!f \|_{W^\mu_2(D)} \le  \underbrace{\|f- I_X\!f \|_{W^\mu_2(D)}}_{A} + \underbrace{\| I_X\!f - \widetilde I_X\!f \|_{W^\mu_2(D)}}_{B}.
\]
Estimating $A$ requires employing $\widetilde D$, $\widetilde X$ and $V_{\widetilde X}$, which were defined in the previous section.  Since $f$ has compact support in $D$, it can be extended to all of $\widetilde D$ by setting it to $0$ in $\widetilde D\setminus D$. Call this extension $f^e$. It follows that $\| f^e \|_{W^\mu_2(\widetilde D \setminus D)} = 0$, so  
\[
\| \tilde f \|_{W^\mu_2(\widetilde D)}=\|f \|_{W^\mu_2(D)}.  
\]
In addition, since $\widetilde X \setminus X \subset \widetilde D\setminus D$, we have that $f^e|_{\widetilde X \setminus X}=0$. Thus, 
\[
I_{\widetilde X} f^e = \sum_{\xi\in \widetilde X} f^e(\xi) \chi_\xi = \sum_{\xi\in X} f(\xi) \chi_\xi + \sum_{\xi\in \widetilde X \setminus X} \underbrace{f^e(\xi)}_0 \chi_\xi = I_X f .
\]
Now, it follows that $\|f- I_X f \|_{W^\mu_2(D)} \le \| f^e - I_X f \|_{W^\mu_2(\widetilde D)} \le \| f^e - I_{\widetilde X} f^e \|_{W^\mu_2(\widetilde D)}$. Since $\widetilde X$ is quasi uniformly distributed in $\widetilde D$, Theorem~\ref{approx_rates} applies, and so we have that
\begin{equation}
\label{bnd_A}
A=\|f- I_X f \|_{W^\mu_2(D)} \le \| f^e - I_{\widetilde X} f^e \|_{W^\mu_2(\widetilde D)} \le Ch^{\beta - \mu} \| f^e \|_{W^\beta_2(\widetilde D)}= Ch^{\beta - \mu} \| f\|_{W^\beta_2(D)}
\end{equation}
To estimate $B$, note that $I_X\!f - \widetilde I_X\!f =  \sum_{\xi \in X} f(\xi)(\chi_\xi(x) - b_\xi(x))$. Applying Lemma~\ref{norm_diff_loc_full}, with $a_\xi = f(\xi)$, we have
\begin{equation}
\label{norm_diff_loc_full_bnd_B}
B= \|\sum_{\xi\in X} f(\xi) (b_{\xi} - \chi_{\xi})\|_{W_2^{\mu}(D)}\le C h^{J-\frac{n}{2} }\|f|_X\|_{\ell_2(X)},
\end{equation}
where $J$ is given in \eqref{J_def}.  Combining this inequality with \eqref{full_riesz_ineq} yields
\begin{equation}
\label{bnd_B_1}
B\le C h^{J-\frac{n}{2} } \|f|_X\|_{\ell_2(X)} \le Ch^{J-n} \| I_Xf \|_{L_2(D)}.
\end{equation}
Furthermore, from \eqref{bnd_A}, with any $h<1$,  we have that 
\[
 \| I_Xf \|_{L_2(D)}\le  \|f- I_Xf \|_{L_2(D)}+ \|f \|_{L_2(D)}\le C\| f\|_{W^\beta_2(D)}+\|f \|_{L_2(D)} \le C\| f\|_{W^\beta_2(D)}.
\]
This and \eqref{norm_diff_loc_full_bnd_B} imply that $B \le Ch^{J-n} \| f\|_{W^\beta_2(D)}$. We can choose $K$ in $J = K\nu/2 +2n -4m-1$ so that $J-n = K\nu/2 +n -4m-1 > \beta - \mu$. From this it follows that $B \le Ch^{\beta -\mu} \| f\|_{W^\beta_2(D)}$. Adding $A$ and $B$ then yields \eqref{q-interp_approx_rate}.
\end{proof}

For future reference, we wish to note that these error estimates lead to estimates for the distance of $f$ to $\spn\{b_\xi: \xi \in X\}$.
Since $\d_{L_2(D)}(f,\spn\{b_\xi: \xi \in X\})\le \|f - \widetilde I_Xf\|_{L_2(D)}$, we have, for $f \in W^\beta_2(D)$ having compact support in $D$,
\begin{equation}
\label{dist_betas}
\d_{L_2(D)}(f,\spn\{b_\xi: \xi \in X\}) \le Ch^{\beta}\|f \|_{W^\beta_2(D)}.
\end{equation}

\begin{remark}
\label{improvements}
\em There are two ways in which Theorem~\ref{q-interp_approx} is likely to be able to be improved: better rates and removal of the  requirement for compact support. As we mentioned earlier, for RBF interpolation of sufficiently smooth functions, Schaback \cite[Theorem~5.1]{Schaback-00-1} obtained a rate double that given earlier in Theorem~\ref{approx_rates}.  Hangebroek \cite[Corollary~5.2]{Hangelbroek-2011-1} established a result showing this phenomenon to be true using the local basis $\{b_\xi\}$, for functions in certain Besov spaces. Something similar is certainly true for Sobolev spaces, and will be dealt with in future work. As to the support requirement, we believe that it is an artifact of the method of proof and is unnecessary, in view of the result \cite[Theorem~4.2]{Narcowich-etal-06-3} holding when all of the centers are inside of $D$. Showing this conjecture holds is an open problem. In section~\ref{numerical_results}, we will discuss numerical evidence supporting our conjectures. 
\em 
\end{remark}

We now turn to a quadrature formula for $f\in W^\beta_2(D)$. We will require this formula to be exact for all functions in $\widetilde V_X$. To derive it, let $s=\sum_{\xi\in X} a_\xi b_\xi$ and observe that this requirement implies that $
\int_D s(x)dx = \sum_{\xi\in X} a_\xi \int_D  b_\xi(x) dx$. If we replace $s$ by the quasi-interpolant $ \widetilde I_X\!f$, then we have 
\begin{equation}
\label{q-interp_quad_def}
Q_X(f) := \int_D \widetilde I_X\!f(x)dx = \sum_{\xi\in X} f(\xi) w_\xi, \  \text{where } w_\xi :=\int_D  b_\xi(x) dx. 
\end{equation}

A straightforward application of Theorem~\ref{q-interp_approx} yields the following error estimates for the quadrature formulas.

\begin{proposition}[\!\!{\cite[Lemma~2]{lero:15}}]
\label{q-interp_quad}
Under the conditions of Theorem~\ref{q-interp_approx}, with $\beta \in \RR$, $n/2<\beta \le m$, $\mu=0$, and $f \in W^\beta_2(D)$ having compact support in $D$, we have 
\begin{equation}
\label{q-interp_quad_est}
\bigg| \int_D f (x)dx - Q_X(f)\bigg| \le  C h^\beta \|f \|_{W^\beta_2(D)}.
\end{equation}
\end{proposition}

We close this section by deriving a formula for the weights in the quadrature formula. In the formula $w_\xi :=\int_D  b_\xi(x) dx$, we replace $b_\xi$ by the right side of \eqref{local_lagrange} and integrate; this yields:
\begin{equation}
\label{weights_calc}
w_\xi = \sum_{\xi\in X} \alpha_{\eta,\xi} \underbrace{\int_D \phi(x - \eta)dx }_{J(\eta)} + \sum_{|\gamma| \le m-1} \beta_{\gamma,\xi\,}\underbrace{\int_D x^\gamma dx}_{J_\gamma}.
\end{equation}
It follows that if we can calculate the $J(\eta)$'s and $J_\gamma$ we can obtain the weights from the coefficients in  \eqref{local_lagrange}. When $D$ is a polygonal domain and $\phi$ a thin-plate spline, there is a simple, exact, analytical formula for $J(\eta)$, which we derive in Appendix~\ref{appendix}. Employing this formula greatly reduces the cost of finding the weights.

\section{Coercivity}
\label{coercivitiy}

In the sequel, we will need various coercivity results for the quadratic form \eqref{bilinear-form}. (At this point, we again use $\Omega$, $\Omegal$, and $\Omega_\cali$ as in section~\ref{var-form}.) We begin with the following lemmas.
\begin{lemma}\label{rho_bound}
Let $u\in L_2(\Omegal )$. Suppose that $0\le \varepsilon <1$. If   $\big |\int_{\Omega_{\mathcal{I}}}u(x)dx\big| \le \varepsilon |\Omega_{\mathcal{I}} |^{1/2} \|u\|_{\Omega_{\mathcal{I}}}$, then
\begin{equation}
\label{rho_inequality}
\frac{1}{\sqrt{|\Omegal  |}}\bigg |\int_{\Omegal }u(x)dx\bigg| 
\le \big(\sqrt{1-\varrho}+\sqrt{\varrho}\,\varepsilon \big)\|u\|_{\Omegal },\ \varrho:=\frac{|\Omega_{\mathcal{I}}|}{|\Omegal |}.
\end{equation}
Furthermore,  if $0<t\le 1$ and $\varepsilon \le \frac{(1-t)\sqrt{\varrho}}{1+\sqrt{1-\varrho}}$, then 
\begin{equation}
\label{rho_eps_inequality}
\frac{1}{\sqrt{|\Omegal  |}}\bigg |\int_{\Omegal }u(x)dx\bigg| 
\le \bigg(1 -  \frac{t\varrho}{1+\sqrt{1-\varrho}}\bigg)\|u\|_{\Omegal }.
\end{equation}
\end{lemma}

\begin{proof}
Since $\int_{\Omegal }u(x)dx = \int_{\Omega}u(x)dx+\int_{\Omega_{\mathcal{I}}}u(x)dx$, by Schwarz's inequality, we have that 
\[
\bigg|\int_{\Omegal }u(x)dx \bigg| \le |\Omega |^{1/2} \|u\|_{\Omega} + \varepsilon |\Omega_{\mathcal{I}} |^{1/2} \|u\|_{\Omega_{\mathcal{I}}} \le \big(|\Omega |^{1/2}+ \varepsilon |\Omega_{\mathcal{I}} |^{1/2}\big) \|u\|_{\Omegal }. 
\]
Divide both sides above by $|\Omega_{\mathcal{I}} |^{1/2}$. Note that $|\Omega | = |\Omegal  |- |\Omega_{\mathcal{I}} |$,  so $|\Omega |/| \Omegal | = 1-\varrho$. The resulting inequality is \eqref{rho_inequality}. The second inequality follows from the first, after a little algebra.
\end{proof}

\begin{lemma}\label{lower_norm_bnd}
Let $u\in L_2(\Omegal )$ and $0< t \le 1$. If $\big |\int_{\Omega_{\mathcal{I}}}u(x)dx\big| \le \varepsilon |\Omega_{\mathcal{I}} |^{1/2} \|u\|_{\Omega_{\mathcal{I}}}$, with $\varepsilon  \le \frac{(1-t)\sqrt{\varrho}}{1+\sqrt{1-\varrho}}$, then
\[
\bigg\|u-|\Omegal |^{-1}\int_{\Omegal }u(x)dx \bigg\|_{\Omegal }^2 \ge \frac{t\varrho}{1+\sqrt{1-\varrho}}\|u\|_{\Omegal }^2.
\]
\end{lemma}

\begin{proof}
Note that $|\Omegal |^{-1}\int_{\Omegal }u(x)dx = \langle u, e_0\rangle_{\Omegal } e_0$, where $e_0=|\Omegal |^{-1/2}$ is a constant unit vector in $L_2(\Omegal )$ and $\langle u, e_0\rangle_{\Omegal } e_0$  is the orthogonal  projection of $u$ onto $e_0$. Hence, $\|u-\langle u, e_0\rangle_{\Omegal } e_0\|_{\Omegal }^2= \|u\|_{\Omegal }^2 -     |\langle u, e_0\rangle_{\Omegal }|^2$. Since 
\[
\langle u, e_0\rangle_{\Omegal } = \frac{1}{\sqrt{|\Omegal  |}}\int_{\Omegal }u(x)dx,
\]
we have, by Lemma~\ref{rho_bound}, that 
\[
\begin{aligned}
\|u-\langle u, e_0\rangle_{\Omegal } e_0\|_{\Omegal }^2 & \ge \bigg(1-\bigg(1 -  \frac{t\varrho}{1+\sqrt{1-\varrho}}\bigg)^2\bigg)\|u\|_{\Omegal }^2 \\
& \ge \bigg(1-\bigg(1 -  \frac{t\varrho}{1+\sqrt{1-\varrho}}\bigg)\bigg)\|u\|_{\Omegal }^2=\frac{t\varrho}{1+\sqrt{1-\varrho}}\|u\|_{\Omegal }^2.
\end{aligned}
\]
\end{proof}

We remark that if $\int_{\Omega_{\mathcal{I}}}u(x)dx = 0$, then \eqref{rho_inequality} becomes
\[
\frac{1}{\sqrt{|\Omegal  |}}\bigg |\int_{\Omegal }u(x)dx\bigg| 
\le \big(\sqrt{1-\varrho}\big)\|u\|_{\Omegal },
\]
and thus the lower bound  in Lemma~\ref{lower_norm_bnd} has the form
\[
\bigg\|u-|\Omegal |^{-1}\int_{\Omegal }u(x)dx \bigg\|_{\Omegal }^2 \ge \varrho\|u\|_{\Omegal }^2.
\]

The point of the lemmas proved above is this. Suppose that we have a subspace $\Pi$ of functions in $L_2(\Omega_\cali )$ with the property that $\dist_{L_2(\Omega_\cali)}(1,\Pi)\le \varepsilon |\Omega_{\mathcal{I}} |^{1/2} $. If we consider all $u\in L_2(\Omegal )$ such that $u|_{\Omega_{\mathcal{I}}}$ is orthogonal to $\Pi$ in $L_2(\Omega_{\mathcal{I}})$, then we have that, for every $p\in \Pi$,
\[
\bigg|\int_{\Omega_{\mathcal{I}}} u dx\bigg| = \bigg| \int_{\Omega_{\mathcal{I}}} u (1-p)dx\bigg| \le \|u\|_{\Omega_{\mathcal{I}}} \|1-p\|_{\Omega_{\mathcal{I}}}.
\]
If we minimize over all $p\in \Pi$, then
\begin{equation}
\label{omega_I_bnd}
\bigg|\int_{\Omega_{\mathcal{I}}} u dx\bigg| \le \|u\|_{\Omega_{\mathcal{I}}} \dist_{L_2(\Omega_{\mathcal{I}})}(1,\Pi)\le \varepsilon |\Omega_{\mathcal{I}} |^{1/2} \|u\|_{\Omega_{\mathcal{I}}}.
\end{equation}
We are now in a position to prove the lower bound for the quadratic form $a(u,u)$. 

\begin{theorem}
\label{lower_bnd_a}
Suppose that $\dist_{L_2(\Omega_{\mathcal{I}})}(1,\Pi)\le \varepsilon |\Omega_{\mathcal{I}} |^{1/2}$ and that, for some $0<t \le 1$, $\varepsilon \le \frac{(1-t)\sqrt{\varrho}}{1+\sqrt{1-\varrho}}$. If $\int_{\Omega_{\mathcal{I}}} u(x) p(x)dx=0$ for all $p\in \Pi$, then
\begin{equation}
\label{lower_est_a}
a(u,u)\ge \frac{t \varrho \underline{\lambda}\delta^{d+2}}{1+\sqrt{1-\varrho}}\|u\|_{\Omegal }^2,
\end{equation}
where $\delta$ and $\underline{\lambda}$ are defined in \cite[Corollary~3.4]{Aksoylu-Mengesha-2010}. 
\end{theorem}

\begin{proof}
By \cite[Corollary~3.4]{Aksoylu-Mengesha-2010} we have that, for all $w\in L_2(\Omegal )$ such that $\int_{\Omegal } w dx = 0$, $a(w,w)\ge \underline{\lambda}\delta^{d+2}\|w\|_{\Omegal }^2$. If $u\in L_2(\Omegal )$, then we have $w=u-|\Omegal  |^{-1}\int_{\Omegal }u(x)dx$ satisfies $\int_{\Omegal }wdx=0$. Furthermore, for any constant $c$, we also have that $a(u-c,u-c)=a(u,u)$. From these facts, we thus have
\[
a(u,u)\ge \underline{\lambda}\delta^{d+2}\big\|u-|\Omegal  |^{-1}\int_{\Omegal }u(x)dx\big\|_{\Omegal }^2.
\]
The lower bound in \eqref{lower_est_a} follows immediately from this inequality, Lemma~\ref{lower_norm_bnd}, and \eqref{omega_I_bnd}.
\end{proof}

\begin{corollary}
\label{lower_est_neumann}
If $\int_{\Omega_\cali}udx =0$, then $a(u,u)\ge \frac{\varrho \underline{\lambda}\delta^{d+2}}{1+\sqrt{1-\varrho}}\|u\|_{\Omegal }^2$. 
\end{corollary}

\begin{proof}
Since $\int_{\Omega_\cali}udx =0$, Lemma~\ref{lower_norm_bnd} applies with $\varepsilon=0$ and $t=1$. The result then follows from the same argument used to prove Theorem~\ref{lower_bnd_a}.
\end{proof}

\section{Lagrange Multiplier Formulation}
\label{lag_mult_form}

We now want to discuss a family of variational problems that will include \eqref{energy} and its discretizations, and these problems into Lagrange-multiplier form. We will deal with the following spaces: $U \subset L_2(\Omegal )$, $\Lambda  \subset L_2(\Omega_{\mathcal{I}})$, and $U^c:= \{u\in U\colon u|_{\Omega_{\mathcal{I}}}\in \Lambda^\perp\}$. All of these are assumed to be closed. 

We also assume that $\Lambda$ satisfies these properties:  First, let $\varepsilon$ satisfy the conditions in Theorem~\ref{lower_bnd_a}. Then, we require that
\begin{equation}
\label{dist_bound}
\dist_{L_2(\Omega_{\mathcal{I}})}(1,\Lambda) \le \varepsilon |\Omega_{\mathcal{I}} |^{1/2}, 
\end{equation}
Second, for every $\lambda\in \Lambda$ there exists an extension\footnote{It might be thought that one can obtain $\tilde \lambda$ by simply taking $\tilde \lambda=0$ on $\Omega$. But since we require $\tilde \lambda\in U$, this will not work in general; however, it will work if $U=L_2(\Omega_\cali)$. See section~\ref{continuous_case}.} $\tilde \lambda\in U$ and a constant $\beta>0$ such that for all $\lambda\in \Lambda$ we have
\begin{equation}
\label{norm_bound}
\beta \le \frac{\|\lambda\|_{\Omega_{\mathcal{I}}}}{\|\tilde \lambda\|_{\Omegal }}.
\end{equation}

Our goal is to establish the following result, which encompasses the various Lagrange multiplier problems that we wish to study.

\begin{proposition} \label{lagrange_form_prop}
There exist unique functions $u\in U$ and $\lambda \in \Lambda$ such that for all $v\in U$ and $\nu\in \Lambda$
\begin{equation}
\label{lagrange_form}
\left\{
\begin{aligned}
a(u,v)+ \int_{\Omega_{\mathcal{I}}} \lambda(x)v(x)dx &= \int_{\Omega} v(x)f(x)dx, \\
\int_{\Omega_{\mathcal{I}}} u(x)\nu(x)dx &=0.
\end{aligned}\right.
\end{equation}
\end{proposition}

The proof will be carried out in several steps. We will begin with the following  inf-sup condition for the linear functional
\[
b(v,\lambda):=\int_{\Omega_{\mathcal{I}}}v(x)\lambda(x)dx,\ v\in U\ \text{and}\ \lambda\in \Lambda.
\]

\begin{lemma} There exists a constant $\beta>0$ such that
\begin{equation}
\label{inf_sup_cond}
\beta \le \inf_{\lambda\in \Lambda} \sup_{v\in U}  \frac{|b(v,\lambda)|}{ \|v\|_{\Omegal }\|\lambda\|_{\Omega_{\mathcal{I}}}}.
\end{equation}
\end{lemma}

\begin{proof}
Let $\lambda$ be fixed. By the second assumption on $\Lambda$,  $\lambda$ has an extension $\tilde \lambda$ to $U$. Because $\tilde \lambda \in U$, we see that
\[
\sup_{v\in U}  \frac{|b(v,\lambda)|}{ \|v\|_{\Omegal }\|\lambda\|_{\Omega_{\mathcal{I}}}} \ge
 \frac{|b(\tilde \lambda ,\lambda)|}{ \|\tilde \lambda\|_{\Omegal }\|\lambda \|_{\Omega_{\mathcal{I}}}}=\frac{\|\lambda \|_{\Omega_{\mathcal{I}}}^2}{ \|\tilde \lambda \|_{\Omegal }\|\lambda \|_{\Omega_{\mathcal{I}}}} = 
 \frac{\|\lambda \|_{\Omega_{\mathcal{I}}}}{ \|\tilde \lambda \|_{\Omegal }}
 \ge \beta,
\]
where the right-hand inequality follows from \eqref{norm_bound}. Taking the infimum above yields \eqref{inf_sup_cond}.
\end{proof}

\begin{lemma} \label{lax-milg_soln}There exists a unique $u_0 \in U^c$ such that $a(u_0,z)=\int_{\Omega}z(x)f(x)dx$ for all $z\in U^c$. 
\end{lemma}

\begin{proof} By Theorem~\ref{lower_bnd_a} and the condition  \eqref{dist_bound}, the quadratic form $a(u,z)$ is coercive on $U^c$; consequently, the Lax-Milgram theorem  implies that $u_0\in U^c$ exists and is unique.
\end{proof}

\begin{proof}[Proof of Proposition~\ref{lagrange_form_prop}]
With $u_0$ in hand, the linear functional below
\begin{equation}
\label{F_functional_def}
F(v):=\int_{\Omega} v(x)f(x)dx - a(u_0,v), \ v\in U,
\end{equation}
satisfies $F(z) = \int_{\Omega}z(x)f(x)dx - a(u_0,z)=0$ for all $z\in U^c$.  In addition, the functional $b$ satisfies satisfies the inf-sup condition \eqref{inf_sup_cond} and is bounded on $U\otimes \Lambda$. By Lemma~10.2.12 in Brenner \& Scott \cite{Brenner-Scott-94-1}, there exists a unique $\lambda_0 \in \Lambda$ such that $b(v,\lambda_0)=F(v)$, where $F$ is given in \eqref{F_functional_def}; that is, 
\[
\int_{\Omega_{\mathcal{I}}}v(x)\lambda_0(x)=\int_\Omega v(x)f(x)dx - a(u_0,v), \ \forall\ v\in U,
\]
so the first equation in \eqref{lagrange_form} holds.  The second is a consequence $u_0$ being in $U^c$. Making the replacements  $u_c \to u$ and $\lambda_0\to \lambda$ completes the proof.
\end{proof}

\subsection{The continuous case with Dirichlet volume constraint}
\label{continuous_case}

We now turn to the problem of solving $a(u,v)=\int_{\Omega}v(x)f(x)dx$, with $u,v=0$ a.e. on $\Omega_{\mathcal{I}}$. Consider the following spaces: $U=L_2(\Omegal )$, $\Lambda =L_2(\Omega_{\mathcal{I}})$, and $U^c := L^c_2(\Omegal) =\{u\in L_2(\Omegal) \colon u |_{\Omega_\cali}=0 \ a.e.\}$. We want to cast this problem into the form \eqref{lagrange_form}.

\begin{theorem}  \label{lagrange_dirichlet_solution}
Let $U$, $\Lambda$ and $U^c$ be as above. Then there exist unique functions $u\in U^c$ and $\lambda\in \Lambda$ that solve \eqref{lagrange_form}.
\end{theorem}

\begin{proof}
We begin by noting that $a(u,v)$ is coercive on $U^c$. This follows from Corollary~\ref{lower_est_neumann}, since all functions in $U^c$ are $0$ on $\cali$, they are obviously orthogonal to $\spn\{1\}$ on $\Omega_\cali$. Moreover, if $\lambda\in \Lambda=L_2(\Omega_\cali)$, then we may define its extension to  $U=L_2(\Omegal )$ by simply setting $\tilde \lambda |_{\Omega} = 0$.  Hence, $\|\tilde \lambda \|_{\Omegal }=\|\lambda\|_{\Omega_\cali}$, and $\Lambda$ satisfies the condition \eqref{norm_bound}, with $\beta=1$. Finally, the condition \eqref{dist_bound} is satisfied, since  $1|_{\Omega_\cali} \in L_2(\Omega_\cali)$ implies that \eqref{dist_bound} holds with $\varepsilon =0$.
\end{proof}

There is an integral-equation approach to this problem. Let $\lambda$ and $u$ be the solutions to the Lagrange equations found above. We start by showing that $\lambda$ is given by an integral operator applied to $u$, and then use this fact to obtain a Fredholm equation for $u$. The assertion concerning $\lambda$ is proved below.
\begin{lemma}
\label{p_formula}
If $\int_{\Omega_{\mathcal{I}}} \lambda(x)\nu(x)dx = - a(u, \tilde \nu), \ \forall \, \nu \in L_2(\Omega_\cali)$, then we have that $\lambda(x) =  -2\int_{\Omega} \gamma(x,y)u(y)dy$, $x\in \Omega_\cali$.
\end{lemma}

\begin{proof} The support of $u$ is $\Omega$. Because $\tilde \nu=0$ on $\Omega$, its support is $\Omega_\cali$. Thus $u(x)\tilde \nu(x)=0$ for all $x\in \Omegal $. This and the symmetry of $\gamma$ then imply that $a(u,\tilde \nu) = -2\int_{\Omegal }\int_{\Omegal } \gamma(x,y) \tilde \nu(x)u(y)dydx$. Using the supports of $u$ and $\tilde \nu$ in the previous expression yields
\begin{equation}
\label{a_u_nu_formula}
a(u,\tilde \nu) = \int_{\Omega_\cali}\bigg( -2\int_{\Omega} \gamma(x,y)u(y)dy\bigg)\nu(x)dx,
\end{equation}
since $\tilde \nu|_{\Omega_\cali}=\nu$. From the definition of $F$ in \eqref{F_functional_def} and $u$ being in $L_2^c(\Omegal)$, we have 
\[
\int_{\Omega_{\mathcal{I}}}\nu(x)\lambda(x)dx = \int_{\Omega_{\mathcal{I}}} \big(-2\int_{\Omega} \gamma(x,y)u(y)dy\big)\nu(x)dx,
\] 
which holds for all $\nu\in L_2(\Omega_\cali)$. Comparing the two sides yields the desired formula for $\lambda$.
\end{proof}

Silling \cite[p.98, eq.~37]{Silling-2010} derives a Fredholm equation of the second kind for a generalization of the type of equilibrium problem we are dealing with here. In our case, the integral equation is the following:
\begin{equation}
\label{peri_fredholm}
\sigma(x)u(x) - \int_{\Omega} \gamma(x,y)u(y)dy = \frac12 f(x), \ \sigma(x)=\int_{\Omegal }\gamma(x,y)dy, \ x\in \Omega.
\end{equation}
The next two results discuss this equation. We begin with the properties of $\sigma$.

\begin{lemma} \label{sigma_lower_bnd}
Let $B_\delta := \{(x,y)\in \Omegal \times \Omegal  \colon |x-y|\le \delta\}$, $\delta>0$. Suppose that there are positive constants $\delta,c_0$ for which $c_0 \le \gamma(x,y)$ for all $(x,y)\in B_\delta$. Then, $\sigma(x)\ge c_0\omega_{n-1} \delta^n/n$, where $\omega_{n-1}$ is the volume of $\sph^{n-1}$.
\end{lemma}

\begin{proof} We may assume that $\delta<\dist(\Omega,\Omegal ^\complement)$. For fixed $x\in \Omega$, the ball centered at $x$ with radius $\delta$ will be in $B_\delta$. Hence, again for fixed $x\in \Omega$, $\gamma(x,y)\ge c_0$, and so $\sigma(x) \ge c_0 \int_{|x-y|\le \delta} dy=\omega_n c_0\delta^n/n$.
\end{proof}

This lemma allows us to divide by $\sigma$, take its square root, and so on. Carrying out such manipulations allows us to put the Fredholm equation \eqref{peri_fredholm} in form given below.

\begin{proposition} With the assumptions made in Lemma~\ref{sigma_lower_bnd}, we have
\begin{equation}
\label{peri_fredholm_explicit}
u(x) - \int_{\Omega} \frac{\gamma(x,y)}{\sigma(x)}u(y)dy = \frac{f(x)}{2\sigma(x)}, \ \ x\in \Omega.
\end{equation}
In addition, if we let $w(x):=\sqrt{\sigma(x)}u(x)$ and $\tilde\gamma(x,y)  = \gamma(x,y)/\sqrt{\sigma(x)\sigma(y)}$, then \eqref{peri_fredholm} has the self-adjoint form
\begin{equation}
\label{peri_sa_fredholm}
w(x) - \int_{\Omega}\tilde \gamma(x,y)w(y)dy=\frac{f(x)}{\sqrt{2\sigma(x)}}, \ x\in \Omega.
\end{equation}
\end{proposition}

For future reference, we point out that when $\gamma(x,y)= \gamma(|x-y|)$ the function $\sigma(x)$ will be constant in $\Omega$. To see this, suppose that the support of $\gamma(r)$ is $[0,\delta]$, where we assume that $\delta<\dist(\Omega,\Omegal ^\complement)$. Fix $x\in \Omega$, the ball $|x-y|\le \delta$ is then contained in $\Omegal $. Thus,
\[
\sigma(x)=\int_{\Omegal } \gamma(|x-y|)dy = \int_{|x-y|\le \delta} \gamma(|x-y|)dy = \omega_{n-1}\int_0^\delta \gamma(r)r^{n-1}dr :=\sigma_\gamma.
\]
The right side is independent of $x$, so $\sigma(x)=\sigma_\gamma$ is constant on $\Omega$.  In fact, it is constant for all $x\in \Omegal $ for which the ball $|x-y|\le \delta$ is contained in $\Omegal $.

\subsection{The discrete case}
\label{discrete_case}
Discretizing the problem begins with choosing a basis of functions to work with. For us, this will be a set of local Lagrange functions associated with a positive definite or conditionally positive definite RBF kernel and a set of centers\footnote{To construct the Lagrange functions, we will make use of a slightly larger set of centers, $Y\supset X$. The centers in $Y\setminus X$ will be outside of $\Omegal$.  } $X$ in $\Omegal $. We will denote the basis by $B=\{b_\xi:\xi\in X\}$. We will use $B$ to construct the three spaces $U_h$, $U^c_h$ and $\Lambda_h$.  As usual, $h$ refers to a mesh norm. We assume that, on $\Omega$, $\Omega_{\mathcal{I}}$, and $\Omegal $, the distribution of centers is quasi uniform. These three spaces are defined this way:  $U_h=\spn\{b_\xi:\xi\in X\}$, $\Lambda_h:=\spn\{b_\xi|_{\Omega_{\mathcal{I}}}: \xi\in X\cap\Omega_{\mathcal{I}}\}$, and $U^c_h=\{u\in U_h\colon u|_{\Omega_\cali} \in \Lambda_h^\perp\}$. 

We now need to discuss conditions \eqref{dist_bound} and \eqref{norm_bound} in connection with $\Lambda_h$. Because RBFs have excellent approximation power, getting $\dist_{L_2(\Omega_\cali} (1|_{\Omega_\cali},\Lambda_h)$ to satisfy the bound in Theorem~\ref{lower_bnd_a} for any $\varepsilon$ only requires taking $h$ sufficiently small and the $K$ in $r_h = Kh |\log(h)|$, sufficiently large. Our next result proves this, along with a coercivity result.

\begin{lemma}\label{RBF_dist_est}
Let $\varrho$, $\varepsilon$ and $t$ be as in Theorem~\ref{lower_bnd_a} and let $\Lambda_h:=\spn\{b_\xi|_{\Omega_{\mathcal{I}}}: \xi\in X\cap\Omega_{\mathcal{I}}\}$ be as in \eqref{bases}, with $\psi_k \to b_\xi$ . Then, for $h$ sufficiently small and $K$ sufficiently large, we have that
\begin{equation}\label{RBF_dist_ineq}
\dist_{L_2(\Omega_\cali)} (1|_{\Omega_\cali},\Lambda_h) \le \varepsilon |\Omega_{\mathcal{I}} |^{1/2}.
\end{equation}
In addition, if $u_h\in U^c_h$, then 
\begin{equation}
\label{discrete_lower_bnd}
a(u_h,u_h)\ge \frac{t \varrho \underline{\lambda}\delta^{d+2}}{1+\sqrt{1-\varrho}}\|u_h\|_{\Omegal }^2.
\end{equation}
\end{lemma}

\begin{proof}
Choose $\alpha>0$ so that the set $\Omega^\varepsilon _{\mathcal{I}}:=\{x\in \Omega_{\mathcal{I}} \colon \dist(x, \partial \Omega_{\mathcal{I}}) \le \alpha \varepsilon\}$ has volume $|\Omega^\varepsilon _{\mathcal{I}}| \le \frac14 \varepsilon^2 |\Omega_{\mathcal{I}} |$. Let $\psi_\epsilon: \Omega_{\mathcal{I}}\to [0,1]$ be a compactly supported $C^\infty$ function  for which $\psi_\varepsilon = 1$ on the set $\Omega_\cali \setminus \Omega^\varepsilon _{\mathcal{I}}$. Next, form the quasi-interpolant $s_h := {\widetilde I}_{X \cap \Omega_{\mathcal{I}}} \psi_\varepsilon \in V_h$.  Applying Theorem~\ref{q-interp_approx}, we have that 
\[
\|\psi_\varepsilon - s_h\|_{L_2(\Omega_\cali)} \le Ch^2 \|\psi_\varepsilon\|_{W^2_2(\Omega_\cali)}.
\] 
for all $h$ sufficiently  small and $K$ sufficiently large. Since $\varepsilon$ is fixed and $h$ and $K$ are at our disposal, we may also choose them so that
\[
Ch^2 \|\psi_\varepsilon\|_{W^2_2(\Omega_\cali)} \le \frac{\varepsilon}{2} |\Omega_{\mathcal{I}} |^{1/2}.
\]
Finally, we note that $\|1|_{\Omega_\cali} - s_h\|_{L_2(\Omega_\cali)}\le \|1|_{\Omega_\cali} - \psi_\varepsilon \|_{L_2(\Omega_\cali)}+ \|\psi_\varepsilon - s_h\|_{L_2(\Omega_\cali)}$. Because $\psi_\varepsilon =1$ on $\Omega_\cali \setminus \Omega^\varepsilon _{\mathcal{I}}$,   we have that $\|1|_{\Omega_\cali} - \psi_\varepsilon \|_{L_2(\Omega_\cali)}=  \|1|_{\Omega_\cali} - \psi_\varepsilon \|_{L_2(\Omega^\varepsilon_\cali)}\le   \|1|_{\Omega_\cali}\|_{L_2(\Omega^\varepsilon_\cali)}\le \frac{\varepsilon}{2} |\Omega_{\mathcal{I}} |^{1/2}$. Hence, $\|1|_{\Omega_\cali} - s_h\|_{L_2(\Omega_\cali)}\le \frac{\varepsilon}{2} |\Omega_{\mathcal{I}} |^{1/2}+\frac{\varepsilon}{2} |\Omega_{\mathcal{I}} |^{1/2}=\varepsilon|\Omega_{\mathcal{I}} |^{1/2}$. The coercivity result \eqref{discrete_lower_bnd} now follows directly from Theorem~\ref{lower_bnd_a}.
\end{proof}

Note that the lower bound in \eqref{discrete_lower_bnd} is independent of $h$, as long as $h$ is sufficiently small -- i.e., $h\le h_0$. This is very important for the method we will use in approximating $u$ by $u_h$. To proceed further, we also need to show that $\Lambda_h$ satisfies the condition in \eqref{norm_bound}. 

\begin{lemma}
\label{norm_bnd_h} 
Consider $\lambda = \sum_{\xi\in X\cap \Omega_\cali} c_\xi b_\xi |_{\Omega_\cali} \in \Lambda_h$ and let $\tilde \lambda :=\sum_{\xi\in X\cap \Omega_\cali} c_\xi b_\xi$, which is an extension of $\lambda$ to $U_h$. Then, there exist constants $\beta>0$ and $h_0>0$, which are independent of $\lambda$, such that $\beta\| \tilde \lambda\|_{\Omegal } \le \| \lambda\|_{\Omega_{\mathcal{I}}}$ holds for all $h\le h_0$.
\end{lemma}

\begin{proof}
From \cite[Theorem~4.11]{HNRW_15} for $h$ sufficiently small, we have 
\[
\|\tilde \lambda\|_{\Omegal } \le C h^{n/2} \| (c_\xi)_{\xi \in X\cap \Omega_\cali}\|_{\ell_2}.
\] 
We will now make use of \cite[Proposition~4.12]{HNRW_15}. Replace $\Omega$ in the proposition by $\Omega_\cali$, $s$ by $\lambda$, and $q$ by $h/ \rho$. Then, we have that 
\[
c\| (c_\xi)_{\xi \in X\cap \Omega_\cali}\|_{\ell_2} \le \rho^{n/2} h^{-n/2} \|\lambda \|_{\Omega_\cali}.
\]
Let $\beta=(Cc\rho^{n/2})^{-1}$. Combining the inequalities then gives $\beta \|\tilde \lambda\|_{\Omegal } \le \|\lambda\|_{\Omega_\cali}$. 
\end{proof}

\begin{theorem}  \label{lagrange_form_h_prop}
Let $U_h$, $\Lambda_h$, $h_0$ and $U^c_h$ be defined as above. For all $h\le h_0$, there exist unique functions $u_h\in U^c_h$ and $\lambda_h \in \Lambda_h$ such that for all $v_h\in U_h$ and $\nu_h\in \Lambda_h$ the following discretized version of \eqref{variationalform} holds:
\begin{equation}
\label{lagrange_form_h}
\left\{
\begin{aligned}
a(u_h,v_h)+ \int_{\Omega_{\mathcal{I}}} \lambda_h(x)v_h(x)dx &= \int_{\Omega} v_h(x)f(x)dx, \\
\int_{\Omega_{\mathcal{I}}} u_h(x)\nu_h(x)dx &=0.
\end{aligned}\right.
\end{equation}
\end{theorem}

\begin{proof} Putting together Lemma~\ref{RBF_dist_est}, Lemma~\ref{norm_bnd_h} and Proposition~\ref{lagrange_form_prop} yields the result.
\end{proof}

\subsection{Error Estimates}\label{sec:Error_Estimates}

To get error estimates,  we will apply results found in sections 10.3 and 10.5 of Brenner \& Scott, which make the assumption that $V_h\subset V$ and $\Pi_h\subset V$. These results hold here because the local Lagrange basis $B=\{b_\xi:\xi\in X\}$ is in $U=L_2(\Omegal )$; and also, the restrictions of them to $\Omega_{\mathcal{I}}$ are in $\Lambda=L_2(\Omega_{\mathcal{I}})$. We can now obtain error estimates for the case at hand.

\begin{theorem} \label{error_est}
Adopt the notation and assumptions made in sections \ref{continuous_case} and \ref{discrete_case}. Then, for $h\le h_0$, 
\[
\|u-u_h\|_{\Omegal } + \|\lambda-\lambda_h\|_{\Omega_\cali} \le C\big(\dist_{L_2(\Omegal )}(u,U_h)+\dist_{L_2(\Omega_\cali)}(\lambda,\Lambda_h)\big).
\]
\end{theorem}

\begin{proof}
Apply Corollary 10.5.18 in Brenner \& Scott.
\end{proof}

At this point getting rates of convergence will depend on two factors: (1) the smoothness of $u$ and $\lambda$; and, (2) the RBF used in the problem. The discussion concerning the Fredholm approach in section~\ref{continuous_case} provides an approach to finding the smoothness of $u$ and $\lambda$. If that can be done, it will be possible to get rates.

The situations for $u$ and $\lambda$ are different. Since $u|_{\Omega_\cali}=0$, the solution $u$ is compactly supported in $\Omegal$ and thus, by Theorem~\ref{q-interp_approx}, the error rates depend only on the smoothness of the kernel $\gamma(x,y)$ and on the source $f$. If these yield $u\in W_2^\beta(\Omegal)$, then distance estimate in \eqref{dist_betas} implies that 
\begin{equation}
\label{u_dist_est}
\dist_{L_2(\Omegal )}(u,U_h) \le Ch^\beta \|u\|_{W_2^\beta(\Omegal)}.
\end{equation}
It may also be possible that $u$ turns out to be in $W_2^{2m}$, then, in view of Remark~\ref{improvements}, we expect that  the error rates should double -- i.e., $h^{2m}$ rather than $h^m$.  This is born out by the numerical results shown in Fig.~\ref{sine_exp}. The expected rate would be about $h^2$, but the rate we obtained is $h^{3.3}$. (It's lower than $h^4$ because $u$ is not quite in $W_2^4$.)

For $\lambda$, the smoothness is known. From Lemma~\ref{p_formula}, we have that $\lambda(x) = -2 \int_{\Omega}  \gamma(x,y) u(y) dy,\ x\in \Omega_\cali$. This formula obviously holds for all $x\in \Omegal$ and thus extends $\lambda$ to $\Omegal$. Differentiating under the integral sign implies that the extension of $\lambda$ to $\Omegal$ is as smooth as $\gamma(x,y)$. 

We also have information about the support of $\lambda$. Since $\gamma(x,y)=0$ for $|x-y|\ge \delta >0$, the Lagrange multiplier $\lambda(x) =0$ when $\d(x\in \Omega_\cali,\Omega) \ge \delta$. Consequently, $\lambda$ is compactly supported in $\Omegal$.  

Unfortunately, this isn't sufficient to apply Theorem~\ref{q-interp_approx} \emph{as stated}. To be able to do that,  $\lambda$ must be compactly supported \emph{in} $\Omega_\cali$. The reason is that the local Lagrange functions employed use only centers in $X\cap \Omega_\cali$, rather than all of $X$.  Even so, as we conjectured in Remark~\ref{improvements}, we expect the to see rates at least those given in Theorem~\ref{q-interp_approx}  to hold. The numerics again bear this out.


\section{Numerical Results}
\label{numerical_results}

We present numerical results for experiments using the discretization described in section~\ref{var-form} and analyzed in sections~\ref{discrete_case} and \ref{sec:Error_Estimates}. The numerical method requires a pre-processing step for constructing the basis, a step of assembling and solving the linear system that arises from the Galerkin method discussed in sections~\ref{discrete_case}., and then a post-processing step for evaluating the $L^2$ error. We discuss the computational methods we employ for each step. The resulting experiments validate the $L^2$ error estimates derived in section~\ref{loc_lag_quasi_interp}. 

 We consider solving two dimensional versions of the problems discussed in section~\ref{discrete_case}, with a radial kernel $\Phi$ and two different  diffusion coefficients $\kappa$; see sections~\ref{STR:LinearAnisotropy} and \ref{STR:ExponentialAnisotropy}. For each experiment,  we test with zero Dirichlet volume constraints although no noticeable difference occurs in the nonzero Dirichlet volume constraint case. The domain of interest for the experiments is denoted $\Omega \cup \Omega_{\mathcal{I}}$ where $\Omega = (0,1) \times (0,1)$ and  $\Omega_{\mathcal{I}} = [-\frac{1}{4},\frac{5}{4}] \times [-\frac{1}{4}, \frac{5}{4}] \backslash \Omega$ denotes the volume constraint region or interaction domain. MATLAB is used for the experiments and plots. Experimental results presented here are the result of directly using the MATLAB backslash operator, which solves the linear set of equations using a sparse direct method.
 
The local Lagrange functions are constructed with linear combinations of the thin plate spline $r^2 \log(r)$. Each local Lagrange function is constructed using approximately $11 \log{N}^2$ nearest neighbor centers, where $N$ is the total number of centers in $\Omega \cup \Omega_I$. The local Lagrange functions are constructed as discussed in section~\ref{sec:local_lagrange_functions}. 

For each numerical experiment, we choose a kernel $\gamma$, an anisotropy term $\kappa$, and a function $u \in L^2_c(\Omega \cup \Omega_I)$ -- i.e., $u$ satisfies the volume constraint. The  source function $f$ is manufactured by computing $\mathcal{L}u(x_i) = f(x_i)$ for each center $x_i$ where
\begin{align}\mathcal{L}u(x) := \int_{\bar{\Omega}} \big(u(x)-u(y)\big) \big( \kappa(x) + \kappa(y) \big) \Phi(\|x-y\|)\, dy
\label{eq:strong_form}
\end{align}
is the strong form corresponding to the bilinear form  \cref{bilinear-form}. We express the kernel $\gamma$ from  \cref{bilinear-form} as $\gamma(x,y) := \big(\kappa(x) + \kappa(y) \big) \Phi(\|x-y\|)$.

The values of $f(x_i)$ are computed by using tensor products of Gauss-Legendre quadrature nodes to approximate the integral in \eqref{eq:strong_form}.

We study $L^2$ convergence of the discrete solution  by constructing sets of uniformly spaced centers with various mesh norms. Uniformly spaced collections of centers $X_h$ are constructed using grid spacing $h = .04, .02, .014$, and $.006$. The convergence of the discrete solution $u_h$ to the solution $u$ is measured by plotting the $L^2$ norm of the error $\|u_h-u\|_{L^2(\Omega \cup \Omega_I)}$ against the mesh norm $h$. The error is computed by placing leveraging tensor products of Gauss-Legendre quadrature nodes over the grid $\bar{\Omega}$.


\subsection{Linear diffusion coefficient}\label{STR:LinearAnisotropy}

We choose $u$, $\kappa$ and the radial function $\Phi$ to be
\begin{align} \label{linear_anisotropy}
\begin{cases}
u(x_1,x_2) = \big(x_1(1-x_1)\big)^{\frac{3}{2}} \big(x_2(1-x_2)\big)^{\frac{3}{2}} \mathbbm{1}_{\Omega}(x_1,x_2)\\
\kappa (x_1,x_2) = 1+x_1+x_2 \\
 \Phi_\varepsilon(\|x-y\|) = \exp\big(-(1- \varepsilon^{-2}\|x-y\|^2)^{-1}\big)\,,
\end{cases}
\end{align}
and thus $\gamma(x,y) := \big(\kappa(x) + \kappa(y) \big) \Phi(\|x-y\|)$, with $x=(x_1,x_2)$ and $y=(y_1,y_2)$.

\begin{figure}
	\centering
	\includegraphics[width=.65\textwidth]{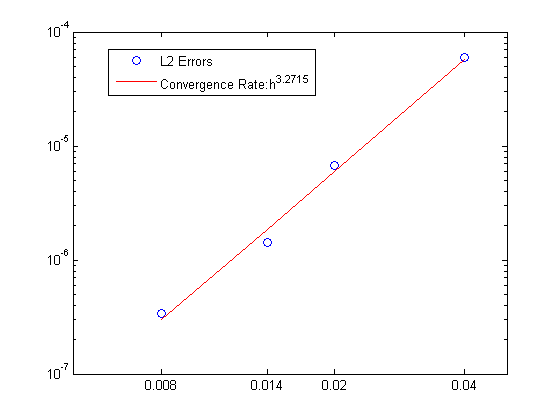}
	\caption{The log of $h$ versus the log of the $L^2$ error for the linear diffusion coefficient experiment with functions given by \eqref{linear_anisotropy} is displayed.}
	\label{polyvar_linear}
\end{figure}

\Cref{polyvar_linear} displays the observed $L^2$ convergence rates with respect to the mesh norm $h$.  The log of the computed $L^2$ error versus the log of the mesh norm is presented along with a best fit line to estimate the convergence order of the observed data.  The smooth solution exhibits a convergence rate of approximately $h^{3}$.

\Cref{TableCondition} displays the condition numbers of the discrete stiffness matrices. The observed condition numbers of the stiffness matrices do not increase as the mesh norm decreases, which suggests that for quasi-uniformly distributed centers, the condition number of the stiffness matrix and the mesh norm $h$ are independent. This prediction is supported by a similar result for the case of a conforming local Lagrange method \cite{lero:15}.

\subsection{Exponential diffusion coefficient
}\label{STR:ExponentialAnisotropy}

For this experiment, we use the functions $u$, $\kappa$ and $\Phi$ given by
\begin{align}\label{exponential_anisotropy}
\begin{cases}
u(x_1,x_2) = \sin(2\pi x_1)\sin(2\pi x_2) \mathbbm{1}_{\Omega}(x_1,x_2) \\
\kappa(x_1,x_2) =\exp(x_1+x_2) \\
 \Phi_\varepsilon(\|x-y\|) = \exp\big(-(1- \varepsilon^{-2}\|x-y\|^2)^{-1}\big)\,.
\end{cases}
\end{align}
Again, $\gamma(x,y) := \big(\kappa(x) + \kappa(y) \big) \Phi(\|x-y\|)$, with $x=(x_1,x_2)$ and $y=(y_1,y_2)$.

\begin{figure}
	\centering
	\includegraphics[width=.65\textwidth]{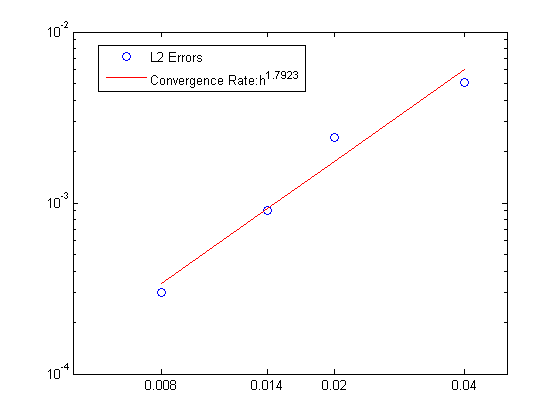}
	\caption{The log of $h$ versus the log of the $L^2$ error for the exponential diffusion coefficient experiment with functions given by \eqref{exponential_anisotropy} is displayed.}
	\label{sine_exp}
\end{figure}

\Cref{sine_exp} displays the $L^2$ convergence plots for the experiments involving $u_2$ and $\kappa_2$. The solution $u$ is not continuously differentiable, so we expect a convergence rate of at most $h^2$. A convergence rate of $h^{1.7}$ is observed.



\begin{table}
\centering
\caption{The mesh norm $h$, number of rows $n$  of the stiffness matrix, and the estimated condition number for the  stiffness matrix with the linear diffusion coefficient \eqref{linear_anisotropy} and the exponential diffusion coefficient \eqref{exponential_anisotropy}. The condition numbers of the stiffness matrices does not increase as $h$ decreases.}
\label{TableCondition}
\begin{tabular}{c c c c}
\multicolumn{2}{c}{} & \multicolumn{2}{c}{Approximate Condition Number} \\ \hline
\noalign{\smallskip}
$h$ &  $n$   & Linear & Exponential  \\ \hline
2.83e-2 & 1444 & 207 & 227 \\
1.41e-2 & 5776 & 171 & 222 \\
9.9e-3 & 11449 & 170 & 219 \\
5.7e-3 & 35344 & 179 & 223 \\
\hline
\end{tabular}
\end{table}

\section*{\bf Acknowledgements}

The authors acknowledge the Texas A\&M University Brazos HPC cluster \cite{brazos_cluster}, which contributed to the research reported here.

\appendix
\section{Appendix}\label{appendix}
In this section we will compute the integrals for the $J(\xi)$'s defined in \eqref{weights_calc}. We will begin by translating $D$ to $D+\xi$, so that in the new coordinates $\xi$ is at the origin and $J(\xi)$ has the form
\[
J(\xi)=\int_{D+\xi} \phi(|x|)dxdy.
\]
To simplify notation, we will use $D$ rather than $D+\xi$, inserting the later at the end of the calculations. 

Suppose that $\phi(|x|)$  satisfies an equation of the form $\Delta\Phi(|x|) = \phi(|x|)$. For example, when $\phi(r)=\phi_2(r)=r^2\log(r)$, we have  $\Phi(r) = \frac{r^4}{32}(2\log(r) - 1)$, where $r=|x|$. When this happens, we may employ Green's theorem to obtain the following formula. 
\[
J(\xi)= \int_{D} \phi(|x| )dxdy =\int_{D}  \Delta \Phi(|x|) dxdy = \oint_{\partial D} \hat \bfn \cdot \nabla \Phi(|x(s)|) ds,
\]
or equivalently,
\begin{equation}\label{stokes}
J(\xi)= \oint_{\partial D} -\frac{\partial \Phi(|x|)}{\partial y}dx+ \frac{\partial \Phi(|x|)}{\partial x}dy, \ |x| =\sqrt{x^2+y^2}.
\end{equation}
Since
\[
\frac{\partial \Phi(|x|)}{\partial x} = \frac{x}{|x|}\Phi'(|x|) \text{ and } \frac{\partial \Phi(|x|)}{\partial y} = \frac{y}{|x|}\Phi'(|x|),
\]
we have
\begin{equation}
\label{line_integral_J}
J(\xi)= \oint_{\partial D} \frac{\Phi'(|x|)}{|x|}\big(-ydx+ xdy\big). 
\end{equation}
It follows that instead of using a 2D quadrature rule, one can get away with a 1D rule. Even better, in the case where $\phi(r)=r^2\log( r)$ and $D$ is a polygonal domain, these integrals can be computed analytically. 

We begin by observing that
\[
\frac{\Phi'(r)}{r} = \frac{r^2}{16}(2\log(r^2)-1),
\]
consequently,
\[
J(\xi)= \oint_{\partial D}\frac{r^2}{16}(2\log(r^2)-1)\big(-ydx+ xdy\big).
\]
If $D$ is a \emph{polygonal} domain, the boundary $\partial D$ consists of a chain of directed line segments. A typical line segment $L$ starts at $(a,A)$  and ends at $(b,B)$. Let $\bfa:=a\bfi  + b\bfj $, $\bfb = b\bfi  + B\bfj$ and $\bfdel:=\bfb - \bfa$. Parametrize $L$ by $\bfx =\bfa+t\bfdel$, $0\le t \le 1$. It is easy to show that $-ydx+ xdy = (aB-bA)dt=(\bfk\cdot\bfa \times \bfdel) dt$. In addition, we have that 
\begin{equation}
\label{parameters}
r^2=|\bfa+t \bfdel|^2 =  \alpha^2+ z^2,\ \text{where } z = |\bfdel| t+ \frac{\bfa\!\cdot\! \bfdel}{|\bfdel|}, \ \alpha:=\pm \sqrt{|\bfa|^2 - \frac{(\bfa\!\cdot\! \bfdel)^2}{|\bfdel|^2}}= \frac{ \bfa\times \bfdel \!\cdot \! \bfk }{|\bfdel |}
\end{equation}
Thus the line integral over $L$ may be put in the form
\[
\int_L \frac{r^2}{16}(2\log(r^2)-1)\big(-ydx+ xdy\big) = \alpha \int_{\frac{\bfa \cdot \bfdel}{|\bfdel|}}^{\frac{\bfb \cdot \bfdel}{|\bfdel|}} \frac{\alpha^2+z^2}{16}(2\log(\alpha^2+z^2)-1)dz.
\]
It follows that we need to compute two indefinite integrals. First, we have 
\[
\begin{aligned}
\int (z^2+ \alpha ^2)\log(z^2+ \alpha ^2)dz & = \int \log(z^2+\alpha^2) d(\frac{1}{3}(z^3+ 3\alpha^2 z)) \\
%
%
=\frac{1}{9}(z^3+3\alpha^2 z)\big(3\log(z^2&+ \alpha ^2) - 2\big) +\frac{4\alpha^2}{3}\tan^{-1}(z/\alpha)-\frac{2\alpha^2}{3}z.
\end{aligned}
\]
And second, $\int (z^2+\alpha^2)dz = \frac{1}{3}(z^3+3\alpha^2 z)$. Combining this result with the previous integral yields
\[
\begin{split}
\int \frac{\alpha^2+z^2}{16}(2\log(\alpha^2+z^2)-1)& dz =  \\ \frac{\alpha z^3+3\alpha^3 z}{144}\big(6\log(z^2+ \alpha ^2)-7\big) +\frac{\alpha^4}{6} \tan^{-1}(z/\alpha) - &  \frac{\alpha^3}{12}z=: f(z,\alpha).
\end{split}
\] 
Finally, we arrive at the integral over the line segment $L$:
\begin{equation}
\label{integral_segment}
\int_L \frac{r^2}{16}(2\log(r^2)-1)\big(-ydx+ xdy\big) =  \alpha \bigg(f\big(\frac{\bfb\!\cdot\! \bfdel}{|\bfdel|},\alpha\big)-f\big(\frac{\bfa\!\cdot\! \bfdel}{|\bfdel|},\alpha\big) \bigg),
\end{equation}
where $\alpha$ is defined in \eqref{parameters}. 

We can give a geometric interpretation to the parameters involved. Let $\hat \bfdel =\bfdel/|\bfdel|$.  Then $\alpha = \bfa \times \hat \bfdel \cdot \bfk$ is the (signed) area of the parallelogram with sides $\bfa$ and $\hat \bfdel$. The endpoints $\bfa \cdot \hat \bfdel $ and $\bfb \cdot \hat \bfdel$ are, respectively,  projections of $\bfa$ and $\bfb$ onto $\bfdel$. 

Restoring $\xi$ to the problem means replacing $D$ above by $D+\xi$, and $L$ by $L+\xi$. The effect on the integrals is to change $\bfa$ and $\bfb$ to $\bfa+\xi$ and $\bfb+\xi$. Of course, $\delta$ remains the same. There is one more step. To get back to the original problem, namely calculating $J(\xi) = \int_{D} \phi(x-\xi)dx$, observe that in a line segment $L_{\text{orig}}$ starting at $\bfa_{\text{orig}}$ and ending at $\bfb_{\text{org}}$, the endpoints are related to those of $L_{\text{orig}}+\xi$ via $\bfa_{\text{orig}}= \bfa+\xi$ and $\bfb_{\text{orig}}= \bfb+\xi$. Thus, in the equations above one should use
\begin{gather}
\bfa = \bfa_{\text{orig}} - \xi \quad \text{and} \quad \bfb = \bfb_{\text{orig}} - \xi, \nonumber \\
\bfdel = \bfb_{\text{orig}} - \bfa_{\text{orig}},  \nonumber  \\
\alpha = \frac{ (\bfa_{\text{orig}} - \xi) \times \bfdel }{| \bfdel |}. \nonumber 
\end{gather}

We conclude by pointing out that the same argument may be used to compute $J(\xi)$ for any TPS  $\phi_m(r) = r^{2m}\log(r)$, $m\ge 1$. Specifically, it is easy to show that 
\[
\Phi_m(r) := \frac{1}{4(m+1)^3}\big( (m+1)\phi_{m+1}(r) - r^{2m+2}\big)
\]
satisfies $\Delta \Phi_m = \phi_m$. Although more complicated, the same integration-by-parts trick still works and will allow us to evaluate $J(\xi)$ exactly.

\bibliographystyle{plain}
\bibliography{rbf}
\end{document}